\documentstyle[12pt]{article}
\newcommand{\be}{\begin{equation}}
\newcommand{\ee}{\end{equation}}
\newcommand{\bqn}{\begin{eqnarray}}

\newcommand{\eqn}{\end{eqnarray}}

\newcommand{\bd}{\begin{description}}
\newcommand{\ed}{\end{description}}

\newtheorem{stat}{}[section]

\def\bs{\begin{stat}}
\def\es{\end{stat}}
\def\ben{\begin{enumerate}}
\def\een{\end{enumerate}}

\def\bp{\noindent{\bf Proof.}  \ \ \ }
\def\ep{\hfill $\Box$}

\makeatletter
\@addtoreset{equation}{section}

\makeatother

\begin{document}

\begin{center}
{\large {\bf On the Cycle Space of a 3--Connected Graph}}
\\[4ex]
{\large {\bf Alexander Kelmans}}
\\[2ex]
{\bf Rutgers University, New Brunswick, New Jersey}
\\[0.5ex]
{\bf University of Puerto Rico, San Juan, Puerto Rico}
\end{center}

\begin{abstract}
We give a simple proof of Tutte's theorem stating 
that the cycle space of a 3--connected graph is 
generated by the set of non-separating circuits 
of the graph. 
\\[1ex]
\end{abstract}

\section{Introduction}

\indent

We consider undirected graphs with no loops and 
no parallel edges. 
All notions on graphs that are  not defined here 
can be found in \cite{D,V}.

Let $G = (V, E, \psi )$ be a graph, where
$V = V(G)$ is the set of vertices, $E = E(G)$ is the set of edges, and $\psi : E \to (^V_2)$ is the edge-vertex incident function.

If $C$ is a cycle of $G$ then $E(C)$ is called 
a {\em circuit} of $G$.
If $X, Y \subseteq E$, then let $X + Y$ denote 
the symmetric difference of $X$ and $Y$, 
i.e. $X + Y = (X \cup Y) \setminus (X \cap Y)$.
Then $2^E$ forms a vector space over $GF(2)$.
Let ${\cal C}(G)$ denote the set of 
circuits of $G$, and so ${\cal C}(G) \subseteq 2^E$. 
Let ${\cal CS}(G)$ denote the subspace of $2^E$ 
generated by ${\cal C}(G)$.
This subspace is called the {\em cycle space of} $G$.
Obviously $X \in {\cal CS}(G)$ if  and only if every 
vertex $v$ in the subgraph of $G$ induced by $X$ 
has even degree.
In particular, $\emptyset \in {\cal CS}(G)$.
If $Z \subseteq E$, then let $G / Z$ ($G \setminus Z$) denote the graph 
obtained from $G$ by contracting (respectively, deleting) the edges in $Z$.
If  $A$ and $B$ are  subgraphs of $G$,
we write, for simplicity,
$G / A$ instead of $G/ E(A)$,  
$A + B$ instead of $E(A) + E(B)$, and 
$A \in {\cal F}$ instead of $E(A) \in {\cal F}$ 
for ${\cal F} \subseteq 2^E$.

A cycle $C$ (the corresponding circuit $E(C)$) in a 
connected graph $G$ is called {\em separating}
if $G /C$ has more blocks than $G$, and
{\em non-separating}, otherwise.
Let ${\cal NC}(G)$ denote the set of non-separating 
circuits
of $G$, and so ${\cal NC}(G) \subseteq {\cal C}(G)$.

Given two graphs $G$ and $F$ with $E(G) = E(F)$, 
we say that  $G$ is {\em strongly isomorphic to} $F$ if
there is an isomorphism $v : V(G) \to V(F)$ from $G$ to 
$F$ that induces the identity map $\epsilon : E \to E$.
\\

One of the classical Whitney theorems states:
\bs {\em \cite{W}}
\label{WhitneyC}
Let $G$ and $F$ be two graphs such that 
$E(G) = E(F)$ and ${\cal C}(G) = {\cal C}(F)$.
If $G$ is 3--connected and $F$ has no isolated vertices, 
then $G$ is strongly isomorphic to $F$.
\es

A very simple proof of {\bf \ref{WhitneyC}} is given in
\cite{K1,K2}.
\\

In \cite{K1} we proved the following strengthening of 
{\bf \ref{WhitneyC}}.
\bs 
\label{WhitneyNC}
Let $G$ and $F$ be two graphs such that 
$E(G) = E(F)$ and ${\cal NC}(G) = {\cal NC}(F)$.
If $G$ is 3--connected and $F$ has no isolated vertices,
then $G$ is strongly isomorphic to $F$.
\es

In \cite{Ksemi} we gave some other strengthenings 
of the Whitney theorem {\bf \ref{WhitneyC}}.
\\

The following theorem, due to W. Tutte \cite{T}  and, independently, A. Kelmans \cite{K1,K2},
is an important result in the study of the graph cycle 
spaces. 
\bs
\label{NCgenerator}
The set of non-separating circuits of a 3--connected 
graph generates the cycle space of the graph.
\es

The above Theorem is an obvious Corollary of
{\bf \ref{WhitneyNC}}.
On the other hand, {\bf \ref{WhitneyNC}} follows from 
{\bf \ref{WhitneyC}} and {\bf \ref{NCgenerator}}.
\\

In  \cite{K1} we proved the following theorem.
\bs
\label{NSandCoind}
Suppose that $G$ is a 3--connected graph,
$X \subseteq E(G)$ and $G \setminus X$ is a connected graph.
Then there exist two distinct non-separating circuits 
$A$, $B$ in $G$ such that 
$|A \cap X| = 1$ and $|B \cap X| = 1$.
\es

We also gave the following simple 
\\[1ex]
{\bf Proof}~~of {\bf \ref{WhitneyNC}}, and therefore also 
{\bf \ref {NCgenerator}}, using {\bf \ref{NSandCoind}} \cite{K1}.
Let $G$ be a 3-connected graph. It is sufficient 
to show that the set ${\cal K}(G)$ of cocircuits 
(i.e. minimal edge cuts) of $G$ is uniquely 
defined by the set ${\cal NC}(G)$ of 
non-separating circuits of $G$. 
Let ${\cal K}'(G)$ be the set of edge subsets $X$ 
of $G$ such that $X \ne \emptyset $  and
$|X \cap C| \ne 1$ for every $C \in {\cal NC}(G)$. 
Obviously ${\cal K}(G) \subseteq {\cal K}'(G)$.
Let ${\cal K}''(G)$ be the set of members of 
${\cal K}'(G)$ minimal by inclusion. 
By {\bf \ref{NSandCoind}} , 
if $X \in {\cal K}'(G)$, then there is $Y \in {\cal K}(G)$ 
such that $Y \subseteq X$. Since $Y \in {\cal K}(G)$, 
every proper subset of $Y$ is not in ${\cal K}(G)$. 
Therefore 
\\
${\cal K}(G) \subseteq {\cal K}'(G) 
\Rightarrow {\cal K}''(G) = {\cal K}(G)$.   
\ep
\\

There are several other proofs of {\bf \ref {NCgenerator}}
(see, for example, \cite{D,V}).
\\

In this paper we give a new fairly simple  proof of 
{\bf \ref {NCgenerator}}.
\\

The results of this paper were presented at the Moscow Discrete Mathematics Seminar in 1977 (see also \cite{Krrr4-06}).

\section{Proof of {\bf \ref {NCgenerator}}}

\indent

We call a graph {\em topologically  3--connected}, or 
simply {\em top 3--connected}, if it is a subdivision of 
a 3--connected graph. A subdivision of a graph $G$ is called {\em top G}.

A {\em thread} in  $G$ is a path $T$ in $G$ such that 
the degree of every inner vertex of $T$ 
is equal to  two and the degree of every end-vertex of $T$ 
is not equal to two in $G$. 
Obviously if $C$ is a cycle of $G$ and 
$E(C) \cap E(T) \ne \emptyset $, then $T \subseteq C$.
If $T$ is a thread in $G$, we write $G - (T)$ instead of 
$G - (T - End(T))$.

A path $P$  with end-vertices $x$ and $y$ is called 
a {\em path-chord of} a cycle $C$ in $G$ if  
$V(C) \cap V(P) = \{x,y\}$, and 
$E(C) \cap E(P) = \emptyset $.
\\

We need the following known facts.
\bs {\em \cite{K2}}
\label{EarAssembly}
Let $G$ be a top  3--connected graph and $G$  not 
top $K_4$.
Then $G$ has a thread $T$ such that
$G - (T)$ is also a top 3--connected graph.
\es

\bs {\em \cite{K2}}
\label{C,C1,C2}
Let $G$ be a top 3--connected graph, 
$C$ a cycle of $G$, and $T$ a thread of $G$ which 
is a path-chord of $C$, and let $R$, $S$ 
be the cycles of $C \cup T$ distinct from $C$.
If $C$ is a non-separating cycle of $G - (T)$, then $R$ 
and $S$ are non-separating cycles of $G$.
\es 

\bp Let $Q = S - (T)$. Then $G / R$ has a block, 
say $H$, containing $E(Q)$. 
Suppose, on the contrary, that $R \not \in {\cal NC}(G)$, 
i.e. $G / R$ has a block $B$ distinct from $H$. 
Then $B$ is also a block of $G / C$.
Suppose that $E(H) \ne E(Q)$. Let $P$ be a block 
of $G / C$ that meets $E(H) \setminus E(Q)$.  
Then $E(P) \ne E(B)$ and  $E(P) \ne E(T)$, and 
therefore
$C \not \in {\cal NC}(G - (T))$, a contradiction.
Thus $E(H) = E(Q)$. 
Then $Q$ is a thread of $G$ and $Q$ is parallel to $T$. Therefore $G$ is not top 3--connected, a contradiction.
\ep
\\

\bs
\label{Theta} {\em \cite{K1,K2}}
Let $G$ be a 3--connected graph. Then for every edge
$e$ of $G$ there are two non-separating cycles
 $P$ and $Q$ of $G$ such that 
$E(P) \cap E(Q) = e$ and $V(P) \cap V(Q) = \psi (e)$.
\es

{\bf Proof}~~(a sketch).
Since $G$ is top 3--connected, there are two cycles 
$R$ and $S$ such that $R \cap S = T$.
Let ${\cal C}_R$ be the set of cycles $C$ in $G$ 
such that $C \cap R = T$, and so $S \in  {\cal C}_R$.
If $C \in  {\cal C}_R$, then let  $\alpha (C)$ be 
the number of edges of the block of $G / C$ 
containing  $E(R - (T))$.
Let $P$ be a cycle in ${\cal C}_R$ such that
$\alpha (P) = \max \{\alpha (C): C \in  {\cal C}_R \}$.
It is easy to show that $P$ is a non-separating cycle 
of $G$.

Applying the above arguments to $R: = P$ and 
$S: = R$, we find another non-separating cycle 
$Q$ of $G$ such that $P \cap Q = T$.
\ep
\\

Now we are ready to prove the following equivalent 
of {\bf \ref{NCgenerator}}.

\bs
Let $G$ be a top 3--connected graph.
Then ${\cal CS} (G)$ is generated by ${\cal NC}(G)$.
\es
\noindent
{\bf Proof}~~ (uses {\bf \ref{EarAssembly}}, {\bf \ref{C,C1,C2}}, and
{\bf \ref{Theta}}).
We prove our claim by induction on the number $t(G)$ 
of  threads of $G$. If $G$ is top $K_4$, then our claim 
is obviously true. So let $t(G) \ge 7$.
By {\bf \ref{EarAssembly}}, $G$ has a thread $T$ 
such that
$G' = G - (T)$ is top 3--connected. 
By the induction hypothesis, 
${\cal CS} (G')$ is generated by ${\cal NC}(G')$.
Obviously if $Q \in {\cal NC}(G')$ and $T$ is not 
a path-chord of $Q$, then $Q \in {\cal NC}(G)$. 
By {\bf \ref{C,C1,C2}}, if $C \in {\cal NC}(G')$, $T$ is 
a path-chord of $C$, and $R$, $S$ are the cycles 
of $C \cup T$ distinct from $C$, then 
$R, S \in {\cal NC}(G)$.
In this case $C = R + S$.
Therefore every cycle in $G'$ is generated by  
${\cal NC}(G)$.
Now let $A$ be a cycle in $G$ but not in $G'$.
Then $T \subseteq A$.
By {\bf \ref{Theta}}, there are $P, Q  \in {\cal NC}(G)$ 
such that $P\cap Q = T$.
Since $T \subseteq A$ and $T \subseteq P$, clearly
$A + P \in {\cal CS}(G')$, and so $A + P$ is generated 
by ${\cal NC}(G)$. Since $(A + P) + P = A$ and
$P \in {\cal NC}(G)$, clearly $A$ is also generated 
by ${\cal NC}(G)$.
\ep
\\

More information on this topic can be found in the expository paper \cite{K3}.

\end{document}